\documentclass[reqno,12pt]{amsart}

\input cyracc.def
\font \tencyr = wncyr10
\def\cyr{\tencyr\cyracc}

\usepackage{color}
\newtheorem{theorem}{Theorem}[section]
\newtheorem{lemma}[theorem]{Lemma}
\newtheorem{corollary}[theorem]{Corollary}

\theoremstyle{definition}

 \theoremstyle{remark}
\newtheorem{remark}[theorem]{Remark}

\newcommand\bR{\mathbb{R}}

\newcommand\cH{\mathcal{H}}

\newcommand\frA{\mathfrak{A}}

\newcommand{\WO}
{\overset{\scriptscriptstyle0}%
{W}\,\!}

\newcommand{\cHO}{\overset{\,\,\scriptscriptstyle0} 
{\cH}\,\!}

\newcommand{\mysection}[1]{\section{#1}
\setcounter{equation}{0}}

\newcommand\cbrk{\text{$]$\kern-.15em$]$}} 
\newcommand\opar{
\text{\,\raise.2ex\hbox{${\scriptstyle |}$}\kern-.34em$($}} 
\newcommand\cpar{%
\text{$)$\kern-.34em\raise.2ex\hbox{${\scriptstyle |}$}}\,}
\newcommand\obrk{\text{$[$\kern-.15em$[$}}

\begin{document}

\title[Parabolic equations in one space dimension]
{On parabolic equations in one space dimension}
\author[N.V. Krylov]{N.V. Krylov}%
\thanks{The work was partially supported by  
NSF Grant DMS-1160569}
\address{127 Vincent Hall, University of Minnesota, Minneapolis,
  MN, 55455}
\email{krylov@math.umn.edu}


\begin{abstract}
 Several negative results are presented
concerning the solvability in Sobolev classes of the 
Cauchy problem for the inhomogeneous  second-order
uniformly 
parabolic equations without lower order terms
in one space dimension. The main coefficient is assumed to be
a bounded measurable function of $(t,x)$
bounded away from zero.
We also discuss upper and lower estimates of certain kind
on the fundamental solutions of such equations.
\end{abstract}

\maketitle

\mysection{Introduction}

We are going to consider functions $u(t,x)$ of two
variables $t,x\in\bR$. Denote $D=\partial/\partial x$.
When it makes sense we write
$$
u_{t}=\partial_{t}u=\frac{\partial u}{\partial t},\quad
u_{x}=Du=\frac{\partial u}{\partial x},\quad
u_{xx}=D^{2}u=\frac{\partial^{2} u}{(\partial x)^{2}},...
$$
 
Let $p\in(1,\infty)$ and denote by $W^{1,2}_{p}$
the set of functions $u$ given on $\bR^{2}$
such that $u,u_{x},u_{xx},u_{t}\in L_{p}(\bR^{2})$.
The norm in $W^{1,2}_{p}$ is introduced in the usual
way. Set
$$
Q=(0,1)\times\bR
$$
and define $\WO^{1,2}_{p}(Q)$ as the set of restrictions on $Q$
of functions in  $W^{1,2}_{p}$ each of which is identically zero
for $t\leq0$.

 For $0<\beta\leq\alpha<\infty $ 
denote by $\frA(\beta,\alpha)$ the set of Borel
functions $a(t,x)$ on $Q$ such that $\beta\leq a(t,x)\leq \alpha$
for all $(t,x)\in Q$. One of our main objects
of investigation is the equation
\begin{equation}
                                                  \label{3.23.1}
u_{t}=au_{xx}+f
\end{equation}
in $Q$.
For given $a\in\frA(\beta,\alpha)$ and $f\in L_{p}(Q)$
we will look for solutions of this equation
in class $W^{1,2}_{p}(Q)$ satisfying zero
initial condition or, in other terms, for solutions
in $\WO^{1,2}_{p}(Q)$.

For the author the main source of
interest in the solvability question of such
simplest one-dimensional equations was the theory
of multidimensional parabolic equations with
coefficients which are only measurable
in time and one spacial variable and, say, just independent of
all other variables. It turns out that is we knew
that \eqref{3.23.1} is solvable in all $\WO^{1,2}_{p}(Q)$,
then the multidimensional version of this result
would be also available and would lead to much more general
results about equations with partially regular coefficients
and with easier proofs than
those, for instance, in \cite{DKim}
and the reference therein. The careful reader can 
see it by following and sometimes slightly changing
the arguments in those references.

However, it turned out that not for all $p\in(1,\infty)$
and measurable $a$
one can guarantee the solvability of \eqref{3.23.1}.
We show here that at least for $p\not\in[3/2,3]$ there are
equations which are not solvable.

This situation is quite different from what is known
for   two-dimensional {\em elliptic\/}
equations with measurable coefficients. 
The fact that generally they are not solvable if $p$ is not close to 2
is well-known and was first demonstrated by N.N. Uraltseva
in 1967 (see \cite{LU}
or \cite{Kr10}). Many more
examples of impossibility of solving elliptic equations
in divergence and non divergence forms
in two space dimensions can be found in \cite{DKim1}.
In our parabolic case we cannot exclude even a part of the range
of $p\in[3/2,3]$, and the author has no idea what is going
on in this range.

 We also provide similar results for divergence type equations.

The third line of our investigation is constructing
estimates from below and from above for solutions
of the Cauchy problem with $f=0$ and the initial data
that is the indicator of an interval.
In the case of the estimates from above we are able
to present essentially sharp estimate
for $a$ in classes $\frA(\beta,\alpha)$ in the full range
of $0<\beta<\alpha<\infty$. In the case of estimates
from below we were only able to cover the case
that $1\geq\beta/\alpha>c$, where $c$ is a certain number,
$c>0$. One can probably go further down to zero
by considering solutions of the Ornstein-Uhlenbeck
equation \eqref{1.20.1} for $\lambda>1$
when  solutions given in an integral form
can be found in \cite{Pe} or \cite{Kr03}. We leave
trying to do this to the interested reader only conjecturing that
the left inequality in \eqref{3.29.03} holds for any $\gamma>1$
if $a_{1\gamma}>0$ is chosen sufficiently small.

Pathological behavior of fundamental solutions
of the Cauchy problem for parabolic equations
  even in one space dimension with continuous coefficient $a$
  was noticed quite a while ago in \cite{Il},
where the fundamental solution of the Cauchy problem
 blows up at a point, say
$(0,0)$, for any $t>0$. Then in \cite{FK} and \cite{Sa}
independently examples were constructed 
again with continuous  $a$ in which
for any $t>0$ the fundamental solution
(as a generalized function, in fact a measure) was just
singular with respect to Lebesgue measure.
We add to this line of research a new information
about the integrals of fundamental solutions over intervals.

The article is organized as follows.
In Section \ref{section 4.19.1} we present our main results.
Section \ref{section 4.19.2} contains general results
about the Ornstein-Uhlenbeck equation
\eqref{1.20.1} when $\lambda$ is arbitrary.
In Section \ref{section 3.31.1} we restrict our attention
to $\lambda<0$ and then use the obtained results
in Section \ref{section 4.19.4} to construct
an essentially sharp barrier from above for the
solutions of the Cauchy problem with $f=0$
and the indicator function of an interval
as the initial data. This barrier
serves in Section \ref{section 3.19.1}
 as a solution of \eqref{3.23.1}
with $f=0$ in class $\WO^{1,2}_{p}(Q)$
for $p\in(1,3/2)$ and this and a duality argument
ruin the hope to build a solvability theory in $W^{1,2}_{p}$
for non divergence type equations and $p\in(1,3/2)\cup(3,\infty)$.

In Section \ref{section 4.19.5} we deal with divergence
type equations and basically use the same barrier
and the observation that the $x$-derivative
of a solution of \eqref{3.23.1} is a solution
of a divergence type equation. The final Section \ref{section 3.31.2}
contains the estimate from below alluded to above.

\mysection{Main results}
                                           \label{section 4.19.1}

The reader understands that equations with
$a$ of class $\frA(\beta,\alpha)$
can  be easily transformed  
 into equations with $a$
 of class $\frA(1,\alpha/\beta)$
or $\frA(\beta/\alpha,1)$  
 by using dilations or contractions
of the $t$-axis. Therefore, we only consider
these two classes of $a$.

\begin{theorem}
                                              \label{theorem 3.23.1}
Let $p\in(1,3/2)$, then there exists an
$\alpha=\alpha(p)\in(1,\infty)$ and a function $a\in\frA(1,\alpha)$,
such that equation \eqref{3.23.1} with $f\equiv0$
has a nonzero unbounded solution in class $\WO^{1,2}_{p}(Q)$.
Furthermore, $\alpha(p)\to1$ as $p\downarrow1$
and $\alpha(p)\to\infty$ as $p\uparrow3/2$.

\end{theorem}

\begin{remark}
This theorem shows that
no matter how small the discontinuities of $a$ are allowed,
there is an $a$ and $p>1$ perhaps very close to $1$ such that the
first
assertion of the theorem holds.

This theorem also implies that there is no $p\in(1,3/2)$
such that the estimate
\begin{equation}
                                             \label{3.29.6}
\sup_{Q}|u|\leq N\|u_{t}-au_{xx}\|_{L_{p}(Q)}.
\end{equation}
holds for any given $a\in\frA(1,\alpha(p))$ with a constant,
perhaps depending on $a$ but independent of $u\in\WO^{1,2}_{p}(Q)$.
Recall that according to the parabolic Alexandrov estimate,
for any $\alpha\in(1,\infty)$ there is a constant $N$
such that \eqref{3.29.6} with $p=2$ holds for all 
 $a\in\frA(1,\alpha)$ and $u\in\WO^{1,2}_{2}(Q)$.

The author does not know what could be the least value of
$p$ for \eqref{3.29.6} to hold for any $\alpha\in(1,\infty)$,
 $a\in\frA(1,\alpha)$, and $u\in\WO^{1,2}_{p}(Q)$ with $N$
depending only on $\alpha$.
\end{remark}

\begin{theorem}
                                              \label{theorem 3.23.2}
Let $p\in(3,\infty)$, then for $\alpha=\alpha(p/(p-1))$
there exists a function $a\in\frA(1,\alpha)$,
such that equation \eqref{3.23.1} for some $f\in L_{p}(Q)$
does not have solutions in class $\WO^{1,2}_{p}(Q)$.

\end{theorem}

\begin{corollary}
                                              \label{corollary 3.24.1}

For any $p\in(1,3/2)\cup(3,\infty)$,  there exists an
$\alpha>1$ and a function $a\in \frA(1,\alpha)$ such that,
for any $N\in(0,\infty)$,
the estimate
\begin{equation}
                                        \label{3.10.02}
\|u_{xx}\|_{L_{p}(Q)}
\leq N\|u_{t}-(\lambda a+1-\lambda)u_{xx}\|_{L_{p}(Q)}
\end{equation}
fails to  hold for all $u\in \WO^{1,2}_{p}(Q)$
and $\lambda\in[0,1]$, that is fails to hold
on the set $ \WO^{1,2}_{p}(Q)\times[0,1]$.
\end{corollary}

Indeed, otherwise the method of continuity
would prove the unique solvability of \eqref{3.23.1}
in class $\WO^{1,2}_{p}(Q)$
for any $f\in L_{p}(Q)$.

The following few results relate to the divergence type
equations. 
Set $\Lambda=(1-D^{2})^{1/2}$, $H^{n}_{p}(\bR)=
\Lambda^{-n}L_{p}(\bR)$, and $\cHO^{1}_{p}(Q)
=\Lambda\WO^{1,2}_{p}(Q)$. For $f\in L_{p}(Q)$ consider the equation
\begin{equation}
                                                  \label{3.23.3}
u_{t}=(au_{x})_{x}+\Lambda f
\end{equation}
in $Q$. Solutions of this equation
will be looked for  in  $\cHO^{1}_{p}(Q)$.
Since $u=\Lambda w$ for a $w\in \WO^{1,2}_{p}(Q)$
the function $u_{t}$ is $ H^{-1}_{p}(\bR) $-valued
and so are $(au_{x})_{x}$ and $\Lambda f$. Hence,
equation \eqref{3.23.1} has perfect sense, and $u_{t}$
is the strong derivative with respect to $t$ of $u$
as an $H^{-1}_{p}(\bR) $-valued function.

\begin{theorem}
                                              \label{theorem 3.23.11}
Let $p\in(1,3/2)$  and $\alpha=\alpha(p)$. Then there exists  
 a function $a\in\frA(1,\alpha)$,
such that equation \eqref{3.23.3} with $f\equiv0$
has a nonzero solution in class $\cHO^{1}_{p}(Q)$.
Moreover, for this solution the functions
 $u_{t},u_{x}$, and $au_{x}$ are continuously
differentiable functions
of $(t,x)$, so that  equation
\eqref{3.23.3} holds in the classical sense  everywhere in $Q$.

\end{theorem}

\begin{theorem}
                                              \label{theorem 3.23.4} 
Let $p\in(3,\infty)$ and $\alpha=\alpha(p)$. Then there exists
  a function $a\in\frA(1,\alpha)$,
such that equation \eqref{3.23.3} for some $f\in L_{p}(Q)$
does not have solutions in  class $\cHO^{1,2}_{p}(Q)$.

\end{theorem}

Similarly to Corollary \ref{corollary 3.24.1}
we have the following.
\begin{corollary}
                                              \label{corollary 3.24.2}

For any $p\in(1,3/2)\cup(3,\infty)$,  there exists an
$\alpha>1$ and a function $a\in \frA(1,\alpha)$ such that,
for any $N\in(0,\infty)$,
the estimate
$$
\|u_{x}\|_{L_{p}(Q)}
\leq N\|f\|_{L_{p}(Q)}
$$
fails to hold on the set of all couples $(\lambda, u)$, where
  $\lambda\in[0,1]$  and  $u\in \cHO^{1}_{p}(Q)$
is a solution  of 
$$
u_{t}=([\lambda a+1-\lambda]u_{x})_{x}+\Lambda f
$$
with $f\in L_{p}(Q)$.
\end{corollary}

\begin{remark}
The author does not know what is going on 
concerning the above results for  the whole region of values of
$p\in[3/2,3]$.
\end{remark}

In spite of the above negative results, for any $\alpha>1$,
if $p$ is sufficiently close to 2, there is a constant $N$ such that
\eqref{3.10.02} holds for all $u\in \WO^{1,2}_{p}(Q)$,
  $\lambda\in[0,1]$, and $a\in \frA(1,\alpha)$. This fact should be
considered well known (in case $p=2$ it is found in
\cite{FK}). It can be easily retrieved from 
Theorem 2.6 of \cite{DK} by following what is said
in Remark 2.3 there. Even better way is to prove it directly
as follows.

By an elementary Lemma 7 of \cite{Kr70}, if $\delta\in(0,1)$, then for any 
$a\in\frA(\delta,\delta^{-1})$, smooth $u(t,x)$ and $p>1$
\begin{equation}
                                             \label{3.29.1}
|u_{t}-u_{xx}|^{p}\leq(1-\delta^{2}/2)^{p}(u_{t}^{2}+u_{xx}^{2})^{p/2}
+(2/\delta)^{p}|u_{t}-au_{xx}|^{p}.
\end{equation}

One also knows that
$$
\|u\|_{1,2,p}:=
\bigg(\int_{Q}(u_{t}^{2}+u_{xx}^{2})^{p/2}\,dxdt\bigg)^{1/p}
$$
defines an equivalent norm in $\WO^{1,2}_{p}(Q)$. Then observe that
by integrating by parts
or, better yet, using the Fourier transform
 we get that for any  $u\in\WO^{1,2}_{2}(Q)$
$$
\int_{Q}(u_{t}-u_{xx})^{2}\,dxdt=
\|u\|^{2}_{1,2,2}+2\int_{Q}u_{t}u_{xx}\,dxdt
$$
$$
=\|u\|^{2}_{1,2,2}+\int_{\bR}|u_{x}(t,\cdot)|^{2}\,dx
\geq\|u\|^{2}_{1,2,2}.
$$
This implies that the norm of the inverse 
operator to  $\partial_{t}-D^{2}:\WO^{1,2}_{2}(Q)\to L_{2}(Q)$
is less than one. The Riesz's convexity theorem implies that
the norm $N_{p}$
 of the inverse to  $\partial_{t}-D^{2}:\WO^{1,2}_{p}(Q)\to 
L_{p}(Q)$ is a continuous function of $p$ and hence
its product with $(1-\delta^{2}/2) $ is less than $\varepsilon$,
which is strictly less than
one for $p$ sufficiently close to 2 (with sufficiently close
defined by $\delta$). Then owing to \eqref{3.29.1} we conclude that
for any $u\in \WO^{1,2}_{p}(Q)$
$$
\int_{Q}|u_{t}-u_{xx}|^{p}\,dxdt
\leq (1-\delta^{2}/2)^{p}N_{p}^{p}\int_{Q}|u_{t}-u_{xx}|^{p}\,dxdt
$$
$$
+\frac{2^{p}}{\delta^{p}}\int_{Q}|u_{t}-au_{xx}|^{p}\,dxdt,\quad
\int_{Q}|u_{t}-u_{xx}|^{p}\,dxdt
$$
$$
\leq 
\frac{2^{p}}{\delta^{p}(1-\varepsilon^{p})}\int_{Q}|u_{t}-au_{xx}|^{p}
\,dxdt,
$$
$$
\|u\|_{1,2,p}\leq 
\frac{2^{p}N^{p}_{p}}{\delta^{p}(1-\varepsilon^{p})}
\int_{Q}|u_{t}-au_{xx}|^{p}\,dxdt.
$$
The latter is an a priori estimate which allows one 
to prove the unique solvability for equation
\eqref{3.23.1} for $p$ close to 2 by the method of 
continuity.
 
The above arguments had, in particular, the goal to be combined with 
the maximum principle and provide for each 
$a\in\frA(\delta,\delta^{-1})$ a transition kernel 
$P_{a}(t,x,s,\Gamma)$ which
is 

(i) a Borel function on $\{1\geq t\geq s\geq0\}\times
\bR$ for any Borel $\Gamma\subset \bR$,

(ii) a
probability measure with respect to $\Gamma$
whenever $1\geq t\geq s\geq0,x\in\bR$,

(iii) for any $\xi\in C^{\infty}_{0}(\bR)$
and $s\in[0,1)$
the function
$$
u(t,x)=\int_{\bR}\xi(y)\,P_{a}(t,x,s,dy)
$$
is a unique continuous in $[s,1]\times\bR$
solution of equation \eqref{3.23.1}
belonging to $W^{1,2}_{2}((s,1)\times\bR)$ and 
satisfying $u(s,x)=\xi(x)$. The details
of construction of $P_{a}(t,x,s,\Gamma)$
can be found in \cite{FK}.

Our next two results concern estimates
on $P_{a}(t,x,s,\Gamma)$.

\begin{theorem}
                              \label{theorem 3.29.1}
Let $\gamma\in(0,1)$. Then
there are   constants  $\nu\in(0,1)$ and $\beta=\beta(\gamma)
\in(1,\infty)$ depending only on $\gamma$
such that for any $\varepsilon\in(0,1/2)$
\begin{equation}
                                      \label{3.29.3}
\nu^{-1}\varepsilon^{\gamma-1}\geq
\frac{1}{\varepsilon}\sup_{a\in \frA(1,\beta(\gamma))}
P_{a}
(1,0,0,[-\varepsilon,\varepsilon])
\geq \nu\frac{\varepsilon^{\gamma-1}}
{|\ln\varepsilon|^{\gamma/2}}.
\end{equation}
Furthermore, $\beta(\gamma)\to1$ as $\gamma\uparrow1$
and $\beta(\gamma)\to\infty$ as $\gamma\downarrow 0$.
 
\end{theorem}

\begin{theorem}
                              \label{theorem 3.29.2}
Let $\gamma\in(1,2)$. Then
there are   constants  $\nu\in(0,1)$ and $\beta=\beta(\gamma)\in(0,1)$ depending only on $\gamma$
such that for any $\varepsilon\in(0,1/2)$
\begin{equation}
                                      \label{3.29.03}
\nu^{-1}\varepsilon^{\gamma-1}\geq
\frac{1}{\varepsilon}\inf_{a\in \frA(\beta(\gamma),1)}
P_{a}
(1,0,0,[-\varepsilon,\varepsilon])
\geq \nu\frac{\varepsilon^{\gamma-1}}
{|\ln\varepsilon|^{\gamma/2}}.
\end{equation}
 Furthermore, $\beta(\gamma)\to1$ as $\gamma\downarrow1$
and $\beta(\gamma)$ tends to a nonzero limit as $\gamma
\uparrow 2$.
\end{theorem}

\mysection{Auxiliary results}
                                           \label{section 4.19.2}
 
In this section we assume that we are given
  $\lambda\in\bR,-\infty\leq a<b\leq
\infty $ and  a function $u>0$ that is
continuous   on $[a,b]\cap\bR$ and satisfies
\begin{equation}
                                                  \label{1.20.1}
u''-2xu'+2\lambda u=0
\end{equation}
on $(a,b)$.

\begin{lemma}
                                    \label{lemma 1.3.1}
The functions $u(-x)$ and
$$
w(x)=u(x)\int_{a}^{x}\frac{1}{u^{2}(t)}e^{t^{2}}\,dt
$$
are solutions of \eqref{1.20.1} on $(-b,-a)$ and $(a,b)$,
respectively.
\end{lemma}

Proof. Both assertions are consequences of direct calculations.
For instance,
$$
w'=u'\int_{a}^{x}\frac{1}{u^{2}(t)}e^{t^{2}}\,dt
+\frac{1}{u}e^{x^{2}},
$$
$$
w''=u''\int_{a}^{x}\frac{1}{u^{2}(t)}e^{t^{2}}\,dt
+2u'\frac{1}{u^{2}}e^{x^{2}}+u(-2u'\frac{1}{u^{3}}
+2x\frac{1}{u^{2}})e^{x^{2}},
$$
and the assertion about $w$ follows.

\begin{lemma}
                                    \label{lemma 1.3.01}
 The function $u$ is  a solution of \eqref{1.20.1} on $(a,b)$
if and only if the function $w(x)=u(x)e^{-x^{2}}$
is a solution of 
\begin{equation}
                                 \label{1.20.2}
w''+2xw'+2\gamma w=0
\end{equation}
on $(a,b)$, where $\gamma=\lambda+1$.
\end{lemma}

The result follows from simple calculations showing that
$$
u(x)=w(x)e^{x^{2}},\quad u'(x)=[w'(x)+2x w(x)]e^{x^{2}},
$$
$$
 u''(x)=[w''(x)+4xw'(x)+(4x^{2}+2)w(x)]e^{ x^{2}}.
$$

\begin{lemma}
                                   \label{lemma 1.2.1}
Let $w\not\equiv0$ be a  solution of \eqref{1.20.2}
on $\bR$ with $\gamma>0$ or $\gamma<-1$,
and let $x_{0}>0$ be a point at which  
$$
w''=0.
$$
 
Then 

(i) $w''>0 $  on any interval $(x_{0},a)$,
$a>x_{0}$, on which
  $w >0$ and $w''<0 $  on any interval $(a,x_{0})$,
$0<a<x_{0}$, on which
  $w >0$;

(ii) $w''<0 $  on any  interval  
$(x_{0},a)$,
$a>x_{0}$, on which
  $w <0$ and $ w ''>0 $  on any interval $(a,x_{0})$,
$0<a<x_{0}$, on which
  $w <0$.

\end{lemma}

Proof. It suffices to prove (i).
Observe that at any point $x>0$ such that $w(x)>0$ and
$w''(x)=0$ we have $x w'=-\gamma w$ and
$$
w'''=-2w'-2\gamma w'=-2(1+\gamma)w'=2\frac{\gamma(1+\gamma)}{x }
w>0.
$$
This easily implies our assertion and the lemma is proved.

\begin{remark}
                                           \label{remark 1.2.1}
Let $w$ be any solution of \eqref{1.20.2},
then $w''(x_{0})=0$ for $x_{0}\ne0$ may only hold
if $w(x_{0})\ne0$, unless $w\equiv0$.
This follows from the uniqueness theorem for ODEs.

\end{remark}

\begin{lemma}
                                            \label{lemma 1.20.3}
 Assume that $v$ and $w$ are solutions
of \eqref{1.20.2} on $[0,\infty)$ and $c_{0}>0$,
$c_{1}>0$ are such that
$$
v(c_{0})\ne0,\quad w(c_{1})\ne0,\quad v''(c_{0})=w''(c_{1})=0.
$$
Define $a_{1}=c_{0}/c_{1}$ and
$$
\psi(x)=v^{-1}(c_{0})v(x),\quad a(x)=1\quad\text{for}\quad x\in[0,c_{0}],
$$
$$
\psi(x)=w^{-1}(c_{1})w(x/a_{1}),\quad a(x)=a_{1}^{2},\quad\text{for}\quad x\in
( c_{0},
\infty).
$$
Then $\psi$ is
 twice continuously differentiable on $[0,\infty)$,
its second-order derivative is Lipschitz continuous
in a neighborhood of $c_{0}$, $a\psi''$ is continuously
differentiable,  $(a\psi'')'$ is Lipschitz continuous
in a neighborhood of $c_{0}$, and $\psi$ satisfies
\begin{equation}
                                      \label{1.20.4}
a(x)\psi''(x)+2x\psi'(x)+2\gamma\psi(x)=0
\end{equation}
on $[0,\infty)$. 

\end{lemma}

Proof. By assumption \eqref{1.20.4}
is satisfied on $[0,c_{0})$. One easily checks that \eqref{1.20.4}
is also satisfied on $(c_{0},\infty)$. Since 
$\psi (c_{0}+)
=\psi (c_{0}-)=1$ and
$\psi''(c_{0}+)
=\psi''(c_{0}-)=0$ we see that  $\psi' (c_{0}+)
=\psi' (c_{0}-) $, so that $\psi$ is twice continuously 
differentiable. Furthermore, $\psi''$ has finite left and right
derivatives  at $c_{0}$, so that it is Lipschitz continuos
near this point. Our assertions concerning   $a\psi''$
follow from the above and \eqref{1.20.4}.
The lemma is proved.

The following lemma shows the way we are going to use
to construct our operators and functions
while proving our main results.
 
\begin{lemma}
                                            \label{lemma 1.20.4}
Under the assumptions of Lemma \ref{lemma 1.20.3}
introduce
 $$
  \Psi(t,x)=\frac{1}{t^{\gamma/2}}\psi
\bigg(\frac{|x|}{2\sqrt{t}}\bigg).
$$
Then $ \Psi\in C^{1,2}_{loc}((0,\infty)\times\bR)$
and $ \Psi$ satisfies
$$
  \Psi_{t}(t,x)=a(t,x) \Psi_{xx}(t,x),
$$
where $a(t,x)=1$ for $|x|\leq 2c_{0}\sqrt{t}$
and $a(t,x)=a_{1}^{2}$ for $|x|> 2c_{0}\sqrt{t}$.

\end{lemma}

This result follows from Lemma \ref{lemma 1.20.3}
and the fact that for $x\geq0$
$$
 \Psi_{t}(t,x)=\frac{1}{t^{1+\gamma/2}}\bigg(-\frac{\gamma}{2}
 \psi\bigg(
\frac{x}{2\sqrt{t}}\bigg)-\frac{1}{2}\frac{x}{2\sqrt{t}}\psi'
\bigg(
\frac{x}{2\sqrt{t}}\bigg)\bigg)
$$
$$
=\frac{a(t,x)}{4t^{1+\gamma/2}} \psi''\bigg(
\frac{x}{2\sqrt{t}}\bigg) ,
$$

\mysection{General properties of solutions
of \eqref{1.20.1} and \eqref{1.20.2} for $\gamma<1$,
$\lambda<0$}
                                          \label{section 3.31.1}

Here we assume that $\gamma<1$,
so that $\lambda<0$.
 One of solutions of \eqref{1.20.1} is
\begin{equation}
                                                  \label{3.22.1}
\phi (x)=\phi_{\lambda}(x)=
\int_{0}^{\infty} e^{-2xr-r^{2}}  {r^{-1-\lambda}}\,dr.
\end{equation}

That $\phi $ is a solution of \eqref{1.20.1}
follows from the fact that
$$
\phi '(x)=
-2 \int_{0}^{\infty} e^{-2xr-r^{2}}  r^{ -\lambda} 
\,dr ,
$$
$$
\phi ''(x)=4
\int_{0}^{\infty} e^{-2xr-r^{2}}  r^{1- \lambda}
\,dr =-2\int_{0}^{\infty} e^{-2xr }  r^{ - \lambda}
\,de^{-r^{2}} 
$$
$$
=2\int_{0}^{\infty}e^{-2xr-r^{2}}[-2xr^{-\lambda}-\lambda
r^{-1-\lambda}]\,dr.
$$

Observe that  the change of variable $r=xs$ and
sending $x\to\infty$ yields
\begin{equation}
                                                \label{1.2.2}
\phi (x)\sim x^{\lambda}
N_{0 },\quad\phi' (x)\sim
  x^{\lambda-1}N_{1 } 
\quad \phi ''(x)\sim  
 x^{\lambda-2}N_{2 } ,
\end{equation}
where 
$$
N_{0}=N_{0\lambda}=\int_{0}^{\infty}e^{-2r}r^{-1-\lambda}\,dr,
$$
$$
N_{1 }=N_{1\lambda}=-2\int_{0}^{\infty}e^{-2r}r^{-\lambda}\,dr
 ,
$$
$$
N_{2 }=N_{2\lambda}=
4
\int_{0}^{\infty} e^{-2r}  r^{1- \lambda}
 \,dr.
$$

To investigate the behavior of 
$\phi (x)e^{-x^{2}}$
as $x\to-\infty$ notice that $\phi (-x)
e^{-x^{2}}$,
$v (x):=[\phi (x)-\phi (-x)]e^{-x^{2}}$,
and
$$
w (x):=\phi (x)e^{-x^{2}}\int_{0}^{x}
\frac{1}{\phi ^{2}(t)}e^{t^{2}}\,dt
$$
are solutions of \eqref{1.20.2}. In addition
\begin{equation}
                                                   \label{1.22.1}
v (x)[2\phi '(0)]^{-1}=
w (x)\phi (0)
\end{equation}
due to the uniqueness theorem for ODEs. Then the formula
\begin{equation}
                                                   \label{3.19.1}
\phi (-x)e^{-x^{2}}
=\phi (x)e^{-x^{2}}-v (x)
\end{equation}
reduces the investigation of the behavior of
$\phi (x)e^{-x^{2}}$
as $x\to-\infty$ to that of $w (x)$ and 
$\phi (x)e^{-x^{2}}$
as $x\to\infty$.

\begin{lemma}
                                            \label{lemma 1.2.2}
As $x\to\infty$,
\begin{equation}
                                                   \label{1.22.2}
v (x)\sim \frac{\phi'(0)\phi(0)}{N_{0 }}
 	x^{-1-\lambda} ,\quad \phi (-x)
\sim  \frac{\phi'(0)\phi(0)}{N_{0 }}
 	x^{-1-\lambda}e^{ x^{2}}.
\end{equation}

\end{lemma}

Proof.
By  \eqref{1.2.2} and l'Hospital's rule
$$
\lim_{x\to\infty}\frac{w (x)}{
x^{-1-\lambda} }=N_{0 }\lim_{x\to\infty}
\frac{\int_{0}^{x} \frac{1}{\phi ^{2}(t)}
e^{t^{2}}\,dt}{
x^{-1-2\lambda}e^{x^{2}}}
$$
$$
=N_{0 }\lim_{x\to\infty}
\frac{1}{\phi ^{2}(x)[-(1+2\lambda)x^{-2-2\lambda}+
2x^{-2\lambda}]  }=\frac{1}{2N_{0 }},
$$
and the first relation in \eqref{1.22.2}
holds due to \eqref{1.22.1}. The second one follows from 
\eqref{3.19.1}. The lemma is proved.

 \mysection{Proof of Theorem \protect\ref{theorem 3.29.1}} 
                                              \label{section 4.19.4}

We split the proof of Theorem  \ref{theorem 3.29.1} into three parts:
first we prove the estimate from above in \eqref{3.29.3},
 then the one from below, and finally we prove its last assertion.
In this section as in  Theorem \ref{theorem 3.29.1},
$\gamma\in(0,1)$, $\lambda=\gamma-1\in(-1,0)$. The function
  $\phi=\phi_{\lambda}$ is taken from \eqref{3.22.1}.
\begin{lemma}
                                            \label{lemma 1.20.8}
Introduce 
$$
u (x)=\phi (x)+\phi (-x).
$$
 Then
there exists a unique $c_{0}=c_{0\lambda }>0$
such that $(u (x)e^{-x^{2}})''=0$ at $x=c_{ 0}$.
In addition, $(u (x)e^{-x^{2}})''<0$ for $|x|< c_{ 0}$
and $(u (x)e^{-x^{2}})''>0$ for $|x|> c_{ 0}$.
Furthermore,
there exists a unique $c_{1}=c_{ 1\lambda}>0$
such that $(\phi (x)e^{-x^{2}})''=0$ at $x=c_{ 1}$.
In addition, $(\phi (x)e^{-x^{2}})''<0$ for $0<x
< c_{ 1}$
and $(\phi (x)e^{-x^{2}})''>0$ for $x> c_{ 1}$.

\end{lemma}

Proof. Observe that at $x=0$
$$
(u (x)e^{-x^{2}})''=2\phi ''(0)-
4\phi (0)=-4(1+\lambda)\phi (0)<0.
$$
In addition, according to \eqref{1.2.2}
and \eqref{1.22.2},
 $u (x)e^{-x^{2}}\to0$ as $x\to
\infty$ and $u (0)>0$. It follows that
the graph of $u (x)e^{-x^{2}}$ has at least
one inflection point on $(0,\infty)$. We denote by
$c_{ 0}$ the smallest one.
Since  $u >0$, Lemma
\ref{lemma 1.2.1} implies that
$(u (x)e^{-x^{2}})''<0$ 
for $0< x < c_{  0}$,
and $(u (x)e^{-x^{2}})''>0$ 
for $ x > c_{  0}$. In particular,
$c_{0\lambda}$ is a unique positive solution
of $(u(x)e^{-x^{2}})''=0$. By symmetry, 
$(u (x)e^{-x^{2}})''<0$ for $|x|< c_{ 0}$
and $(u (x)e^{-x^{2}})''>0$ for $|x|> c_{ 0}$
The same argument (apart from symmetry) works for 
$\phi (x)e^{-x^{2}}$ and the lemma is proved. 

Introduce
$$
a_{1}=a_{1\gamma}=c_{0 }/c_{1 }=c_{0\lambda}/c_{1\lambda},
$$
$$
\hat a(x)=\hat a_{\gamma}(x)=1\quad\text{for}
\quad |x|\leq c_{0},
$$
$$
\hat a(x)=\hat a_{\gamma}(x)=a_{1\gamma}^{2}=:\beta(\gamma)
\quad\text{for}
\quad |x|> c_{0},
$$
$$
w(x)=w_{\gamma}(x)=u^{-1}(c_{0})
u(x)e^{-x^{2}+c_{0}^{2}}\quad\text{for}
\quad |x|\leq c_{0},
$$
$$
w(x)=w_{\gamma}(x)=\phi ^{-1}(c_{1})
\phi (|x|/a_{1} )e^{-(x/a_{1})^{2}+c_{1}^{2}}
\quad\text{for}
\quad |x|> c_{0}.
$$
$$
\Psi(t,x)=\Psi^{(\gamma)}(t,x)=\frac{1}{t^{\gamma/2}}w(x/(2\sqrt{t})),
$$
\begin{equation}
                                       \label{3.22.3}
a^{*}(t,x)=a^{*}_{\gamma}(t,x)=
\hat a(x/(2\sqrt{t})).
\end{equation}

Obviously,
$$
a^{*}_{\gamma }\in\frA(1,\beta(\gamma)).
$$

\begin{remark}
                                                \label{remark 3.22.1}
Lemma \ref{lemma 1.20.3} implies that $w(x)$ has three bounded
derivatives, which obviously tend to zero exponentially fast as $|x|
\to\infty$. This yields that $\Psi$ has three   derivatives
in $(t,x)$ which are bounded in each set $\{t>\varepsilon\}$,
where $\varepsilon>0$ and tend to zero 
exponentially fast as $|x|\to\infty$
provided that $t$ is restricted to a bounded interval separated
from zero.

 Lemma \ref{lemma 1.20.3}
also implies that $\hat aw''(x)$ has two bounded
derivatives, which obviously tend to zero exponentially fast as $|x|
\to\infty$. This yields that $a^{*}\Psi_{xx}$ has two
  derivatives
in $(t,x)$ which are bounded in each set $\{t>\varepsilon\}$,
where $\varepsilon>0$ and tend to zero exponentially fast
 as $|x|\to\infty$
provided that $t$ is restricted to a bounded interval separated
from zero.

\end{remark}

Observe also that by Lemma \ref{lemma 1.20.4}
\begin{equation}
                                       \label{1.15.2}
\Psi_{t}(t,x)=a^{*}(t,x) \Psi_{xx}(t,x)
\end{equation}
and by Lemma \ref{lemma 1.20.8}
\begin{equation}
                                       \label{1.15.3}
  \Psi_{xx}(t,x)\leq0\quad\text{for}\quad
|x|\leq 2c_{0}\sqrt{t},\quad
  \Psi_{xx}(t,x)\geq0\quad\text{for}\quad
|x|\geq 2c_{0}\sqrt{t}.
\end{equation}

\begin{lemma}
                                \label{lemma 1.15.1}
We have $c_{1}<c_{0}$, $a_{1}>1$, and
\begin{equation}
                                       \label{1.15.1}
\Psi_{t}(t,x)=\max_{a\in[1,a^{2}_{1\gamma}]}[a\Psi_{xx}(t,x)].
\end{equation}

\end{lemma}

Proof. Assume that $c_{1}\geq c_{0}$. Then $a_{1}\leq1$
and  \eqref{1.15.2}, \eqref{1.15.3} imply that
$$
\Psi_{t}(t,x)=\min_{a\in[ a^{2}_{1},1]}[a\Psi_{xx}(t,x)],
$$
$$
\Psi_{t}(t,x)\leq \Psi_{xx}(t,x).
$$
It follows by the maximum principle that
$$
\frac{1}{2\sqrt{\pi}}\int_{\bR}\Psi (t,y)
e^{-(x-y)^{2}/4}\,dy\geq \Psi(t+1,x),
$$
$$
\int_{\bR}\Psi(t,y)\,dy\geq 2\sqrt{\pi}\Psi(t+1,x).
$$
However, the integral on the left equals
$$
2t^{(1-\gamma)/2}\int_{\bR}w(y)\,dy\to0
$$
as $t\downarrow0$. This yields a
contradiction, hence $c_{1}<c_{0}$, and the rest is trivial.
 The lemma is proved.

 We now prove the estimate from above in \eqref{3.29.3}.
Recall that $\beta(\gamma)=a_{1\gamma}^{2}$
 for $0<\gamma<1$.

\begin{theorem}
                              \label{theorem 1.15.2}
For any $a\in \frA(a_{1\gamma}^{2})$ and $\varepsilon\in(0,1)$
$$
\frac{1}{\varepsilon}P_{a}(1,0,0,[-\varepsilon,\varepsilon])
\leq N\varepsilon^{\gamma-1},
$$
where the constant $N$ depends only on $\gamma$.
\end{theorem} 

Proof. Since for any $t_{0}>0$ the function 
$\Psi(t_{0}+t,x)$ satisfies
$$
\Psi_{t}(t_{0}+t,x)\geq a(t,x)\Psi_{xx}(t_{0}+t,x),
$$
by the maximum principle
we have
$$
\int_{\bR}\Psi(t_{0},y)\,P_{a}(1,x,0,y)\leq \Psi(t_{0}+1,x),
$$
$$
t_{0}^{-\gamma/2}\int_{\bR}w( y/(2\sqrt{t_{0}}))
\,P_{a}(1,0,0,y)\leq \Psi(t_{0}+1,0),
$$
$$
t_{0}^{-\gamma/2}\min_{|y|\leq1/2}w(y)
P_{a}(1,0,0,[-\sqrt{t_{0}},\sqrt{t_{0}}])\leq \Psi(t_{0}+1,0),
$$
and this proves the theorem.

Next, we prove the estimate from below in \eqref{3.29.3}.

\begin{theorem}
                              \label{theorem 1.16.1} 
There is a constant  $\nu>0$ depending only on $\gamma$
such that for any $\varepsilon\in(0,1/2)$
\begin{equation}
                                      \label{3.22.6}
\frac{1}{\varepsilon}\sup_{a\in \frA(1,a^{2}_{1\gamma})}
P_{a}
(1,0,0,[-\varepsilon,\varepsilon])
\geq \nu\frac{\varepsilon^{\gamma-1}}
{|\ln\varepsilon|^{\gamma/2}}.
\end{equation}
 
\end{theorem}

Proof. Fix a $t_{0}\in(0,1)$ and set 
$a (t,x)=a^{*}(t_{0}
+t,x)$, where $a^{*}$ is introduced in \eqref{3.22.3}.
 By comparing the equations
satisfied by both sides of the following equation
we come to a proof of the fact that
\begin{equation}
                                      \label{1.16.1}
t_{0}^{-\gamma/2}\int_{\bR}
w( y/(2\sqrt{t_{0}}))\,P_{a }
(t,x,0,dy)= \Psi(t_{0}+t,x).
\end{equation}
 It follows due to  \eqref{1.2.2} that
 for any $c\geq c_{0}$ 
we have
$$
t_{0}^{-\gamma/2}w(0)P_{a}(1,0,0,[-2c\sqrt{t_{0}},
2c\sqrt{t_{0}}])
+Nt_{0}^{-\gamma/2  } c ^{\lambda}
e^{-(c /a_{1})^{2}}\geq \Psi(t_{0}+1,0),
$$
where   
 $N$ is a constant depending only on $\gamma$.
For 
$$
c=a_{1}\sqrt{ (\gamma/2)|\ln t_{0}|}
$$
and $t_{0}$ small enough the second term on the left
is less than $(1/2) \Psi(2,0)$ due to the
 fact that $\lambda<0$. Hence, with a constant $\alpha
\geq1$
depending only on $\gamma$ for all small $t_{0}$
\begin{equation}
                                      \label{3.22.5}
t_{0}^{-\gamma/2}w(0)P_{a}(1,0,0,[-\alpha
\sqrt{t_{0}|\ln t_{0}|},
\alpha 
\sqrt{t_{0}|\ln t_{0}|}])
\geq (1/2) \Psi(2,0).
\end{equation}
Now  denote
$$
\varepsilon=\alpha \sqrt{t_{0}|\ln t_{0}|}.
$$
Then
$$
t_{0}=\frac{\varepsilon^{2}}{{\alpha^{2}|\ln 
t_{0}|}},\quad
\ln\varepsilon=\ln\alpha+(1/2)\ln t_{0}
+(1/2)\ln|\ln t_{0}|,
$$
and $|\ln t_{0}|
\leq|\ln\varepsilon|$ for $t_{0}$ small enough, so that
$$
t_{0}\geq \frac{\varepsilon^{2}}{{\alpha^{2}
|\ln \varepsilon|}}.
$$
This allows us to transform \eqref{3.22.5}
into \eqref{3.22.6} for small $\varepsilon$,
for which it only makes any real sense,
and
proves the theorem.

By comparing the behavior of $\varepsilon^{\gamma-1}$
for different $\gamma$ and using the results of
Theorems \ref{theorem 1.15.2} and \ref{theorem 1.16.1}
we immediately come to  the following.

\begin{corollary}
                                            \label{corollary 4.19.1}
The function $\beta(\gamma)=a^{2}_{1\gamma}$
is an increasing function of $\gamma\in(0,1)$.
\end{corollary}

Finally we deal with the last assertion of
Theorem \ref{theorem 3.29.1}.
It suffices to prove that 
\begin{equation}
                                                   \label{4.4.3}
\lim_{\gamma\uparrow 1}c_{1\gamma}=2^{-1/2},
\quad\lim_{\gamma\uparrow 1}c_{0\gamma}=2^{-1/2},
\end{equation}
\begin{equation}
                                                   \label{4.9.1}
\lim_{\gamma\downarrow 0}c_{0\gamma}/c_{1\gamma}=
\infty.
\end{equation}

Observe that if a function $v>0$ satisfies \eqref{1.20.1}, then
the inequality $(ve^{-x^{2}})''>0$ is written as
$$
v''-4xv'+2v(2x^{2}-1)>0,\quad -2xv'+2v(2x^{2}-1-\lambda)>0,
$$
\begin{equation}
                                                   \label{4.5.1}
\frac{xv'(x)}{v(x)}<2x^{2}-1-\lambda.
\end{equation}

Accordingly, since $\phi$ satisfies \eqref{1.20.1},
 equation $(\phi e^{-x^{2}})''=0$, 
defining a unique $c_{1}=c_{1\gamma}>0$, transforms into  
\begin{equation}
                                                   \label{4.4.2}
\frac{c_{1}\phi'(c_{1})}{\phi(c_{1})}=2c_{1}^{2}-1-\lambda.
\end{equation}
Since $\phi'<0$ we conclude that the right-hand side of
\eqref{4.4.2} is negative and $c_{1\gamma}$ is a bounded 
function of $\gamma\in(0,1)$. Furthermore, obviously
$\phi'$ is bounded on any interval $[0,b]$ by a constant
independent of $\gamma\in[1/2,1)$ and $\phi\to\infty$
uniformly on any such interval as $\gamma\uparrow1$,
so that $\lambda\uparrow0$. Now the first 
relation in \eqref{4.4.3} follows from \eqref{4.4.2}.
 
To prove the second one, observe that $c_{0}=c_{0\gamma}$ is defined 
as a unique positive solution of
\begin{equation}
                                                   \label{4.4.4}
c_{0}\frac{\phi'(c_{0})-\phi'(-c_{0})}
{\phi(c_{0})+\phi(-c_{0})}=2c_{0}^{2}-1-\lambda
\end{equation}
and
$$
c_{0}\phi'(c_{0})-c_{0}\phi'(-c_{0})=4c_{0}\int_{0}^{\infty}
e^{-r^{2}}r^{-\lambda}\sinh(2c_{0}r)\,dr,
$$
$$
\phi(c_{0})+\phi(-c_{0})=2\int_{0}^{\infty}
e^{-r^{2}}r^{-1-\lambda}\cosh(2c_{0}r)\,dr.
$$
Furthermore, from the above argument
concerning 
\eqref{4.5.1} and Lemma \ref{lemma 1.20.8} we know that if
$x>0$ and
$$
x\frac{\phi'(x)-\phi'(-x)}
{\phi(x)+\phi(-x)}<2x^{2}-1-\lambda,
$$
then $c_{0}\leq x$. For $x=1$ the left-hand side obviously
tends to zero
as $\gamma\uparrow1$ ($\lambda\uparrow0$). This yields
the boundedness of $c_{0\gamma}$ and we can finish
the proof of
the second relation in \eqref{4.4.3} as before.

To prove \eqref{4.9.1} we   look at \eqref{4.4.2} 
as quadratic equation relative to $c_{1}$ and then find that
$$
c_{1}=2(1+\lambda)/B,
$$
where (recall that $\phi'<0$)
$$
B=B_{\lambda}=|\phi'(c_{1})|\phi^{-1}(c_{1})
+\sqrt{|\phi'(c_{1})|^{2}\phi^{-2}(c_{1})+8(1+\lambda)}.
$$
Since $|\phi'(x)|\phi^{-1}(x)\not\to0$ as $x\downarrow0$,
we conclude that
$$
c_{1\gamma}= O(\gamma)
$$
as $\gamma\downarrow0$.

Now if we assume that   along a sequence $\gamma_{n}
\downarrow0$ we have $c_{0\gamma_{n}}/c_{1\gamma_{n}}\to\beta<\infty$,
then $c_{0\gamma_{n}}\to0$,
 and after dividing both 
parts of \eqref{4.4.4} by $c_{0\gamma}^{2}$ and passing to the limit
along the subsequence we get that
\begin{equation}
                                             \label{4.10.1}
\frac{\phi''_{-1}(0)}{\phi_{-1}(0)}=2-\lim_{n\to\infty}\frac{\gamma_{n}}{
c^{2}_{0\gamma_{n}}}.
\end{equation}
The function $\phi_{-1}$ satisfies \eqref{1.20.1}
with $\lambda=-1$, therefore, $\phi''_{-1}(0)=2\phi_{-1}(0)$,
and \eqref{4.10.1} implies that
$$
\lim_{n\to\infty}\frac{\gamma_{n}}{
c^{2}_{0\gamma_{n}}}=0,\quad
\lim_{n\to\infty}\frac{c_{1\gamma_{n}}}{
c^{2}_{0\gamma_{n}}}=0,
 \quad
\frac{1}{\beta^{2}}\lim_{n\to\infty}\frac{1}{
c _{1\gamma_{n}}}=0,
$$
and the latter is impossible. This proves \eqref{4.9.1}
and brings the proof of Theorem \ref{theorem 3.29.1}
to an end.

\mysection{Proof of Theorems \protect\ref{theorem 3.23.1} and
\protect\ref{theorem 3.23.2}}
                                          \label{section 3.19.1}

For $p\in(1,3/2)$ one can define a function $\gamma(p)\in(0,1)$
so that 
$$
  \gamma(p)<\frac{3}{p}-2,\quad \lim_{p\downarrow1}\gamma(p)=1.
$$
We take any such function $\gamma(p)$ and set
$$
\alpha(p)=\beta(\gamma(p))=a^{2}_{1\gamma(p)}.
$$
Obviously,
$$
a^{*}_{\gamma(p)}\in\frA(1,\alpha(p)).
$$

Theorem  \ref{theorem 3.23.1}  
is a   consequence 
of the following.

\begin{theorem}
                                 \label{theorem 3.10.1}
Let $p\in(1,3/2)$. Then  
the equation
$$
u_{t}=a^{*}_{\gamma(p)}u_{xx}
$$
has a nonzero unbounded solution of class $\WO^{1,2}_{p}(Q)$.
Furthermore, $\alpha(p)\to1$ as $p\downarrow1$
and $\alpha(p)\to\infty$ as $p\uparrow3/2$

\end{theorem}

Proof. Simple computations show that
if we extend $\Psi^{(\gamma(p))}(t,x)$ as zero
for $t\leq0$, then the unbounded function we obtain
will belong to $W^{1,2}_{p}((-1,1)\times\bR)$
and, to prove the first assertion,
 it only remains to recall that $\Psi^{(\gamma)}$
satisfies \eqref{1.15.2}. The second assertion
follows immediately by the construction of $\gamma(p)$
and Theorem \ref{theorem 3.29.1}. The theorem
is proved.

Next argument is based
on a general rule that if a linear homogeneous equation has   
a nonzero solution $u$, then the adjoint equation
is only solvable if its right-hand side
is orthogonal to $u$. We apply this rule to $u=\Psi_{xx}$.

 Here is a result implying Theorem \ref{theorem 3.23.2}.

\begin{theorem}
                                 \label{theorem 3.10.01}
Let $p\in(3,\infty)$. Then for $q=p/(p-1)$ and $\alpha
=a^{2}_{1\gamma(q)}$ there exists
  a function $a\in \frA(\alpha)$,
and an $f\in L_{p}(Q)$ such that
the equation
\begin{equation}
                                        \label{3.19.4}
u_{t}=au_{xx}+f
\end{equation}
has no solutions of class $\WO^{1,2}_{p}(Q)$.

\end{theorem}

Proof. Define  
$$
\Phi(t,x)=\Psi^{(\gamma(q))}(1-t,x),\quad
a(t,x)=a^{*}_{\gamma(q)}(1-t,x),
$$
and take any $f\in L_{p}(Q)$ such that
$\Phi f\in L_{1}(Q)$ and
$$
\int_{Q}\Phi_{xx} f\,dxdt\ne0.
$$

Assume that $u$ is a solution of \eqref{3.19.4} of
class $\WO^{1,2}_{p}(Q)$. We
multiply  \eqref{3.19.4} through
by $\Phi_{xx}$ and use that $a\Phi_{xx}=-\Phi_{t}$.
Then we have the following
$$
\int_{0}^{t_{0}}\int_{\bR}\Phi_{xx} u_{t}\,dxdt
=-\int_{0}^{t_{0}}\int_{\bR}\Phi_{t}u_{xx}\,dxdt+
\int_{0}^{t_{0}}\int_{\bR}\Phi_{xx} f\,dxdt,
$$
$$
\int_{0}^{t_{0}}\int_{\bR}\frac{\partial}{\partial t}[\Phi_{xx} u]
\,dxdt=\int_{0}^{t_{0}}\int_{\bR}\Phi_{xx} f\,dxdt=:\kappa(t_{0}),
$$
\begin{equation}
                                                     \label{5.4.1}
 \int_{\bR}u_{x}(t_{0},x)\Psi_{x}(1-t_{0},x)\,dx
= \kappa(t_{0}),
\end{equation}
where the last equality is obtained by integrating by parts
and using  Remark \ref{remark 3.22.1}.

Next, we are going to use
an embedding theorem that is Lemma 2.3.3 of \cite{LSU}
according to which for each $t\in[0,1]$,
$u(t,\cdot)\in C^{1+\varepsilon}(\bR)$
and the norm of $u(t,\cdot)$ in this space
is a bounded function of $t$. Here $\varepsilon>0$
is any number such that 
$$
1-\frac{3}{p}>\varepsilon>\gamma(q).
$$
That such an $\varepsilon$ exists follows from the fact that
the  inequalities $1-3/p>\gamma(q)$ and $ \gamma(q)<3/q-2$
are equivalent.

Since $\Psi_{x}(1-t_{0},x)$ is an odd function of $x$
we can replace $u_{x}(t_{0},x)$ in \eqref{5.4.1}
 with $u_{x}(t_{0},x)-u_{x}(t_{0},0)$,
and since the latter by magnitude is less than $N|x|^{\varepsilon}$,
where $N$ is a constant,
we come to the conclusion that
$$
|\kappa(t_{0})|\leq \frac{N}{(1-t_{0})^{(1+\gamma)/2}}
\int_{\bR}|x|^{\varepsilon}|w'(x/(2\sqrt{1-t_{0}}))|\,dx,
$$ 
where $\gamma=\gamma(q)$.

However, the change of variable  $y=x/(2\sqrt{1-t_{0}})$
shows that the expression on the right equals
a constant times $(1-t_{0})^{(\varepsilon-\gamma) /2}$,
which tends to zero as $t_{0}\uparrow 1$. Hence,
$\kappa(t_{0})\to 0$ and this is the desired contradiction
proving the theorem.

\mysection{Proof of Theorems \protect\ref{theorem 3.23.11} 
  and
\protect\ref{theorem 3.23.4} }
                                           \label{section 4.19.5}

  To prove
Theorem  \ref{theorem 3.23.11} it suffices to take 
$$
a=a^{*}_{\gamma(p)},\quad u= \Psi_{x}^{(\gamma(p))}
$$
and use what is said in the proof of  Theorem \ref{theorem 3.10.1} 
and also recall the definition of $\cHO^{1}_{p}(Q)$ and use 
Remark \ref{remark 3.22.1}.

While proving Theorem \ref{theorem 3.23.4} we
define $\Phi$ and $a$ as in the proof of Theorem
\ref{theorem 3.10.01} and take any $f\in L_{p}(Q)$
such that $f\Lambda \Phi_{x}\in L_{1}(Q)$ and
$$
\int_{Q}f\Lambda\Phi_{x}\,dxdt \ne0
$$
and let $u$ be a solution of \eqref{3.23.3} 
of class $\cHO^{1,2}_{p}(Q)$.

Then observe that $\Phi_{x}(t,\cdot)$ is
a strongly
continuous differentiable $H^{1}_{q}(\bR)$-valued function  
on $[0,1)$ for any $q\in(1,\infty)$, in particular,
for $q=p/(p-1)$. Since $u(t,\cdot)$ 
is
a strongly
  differentiable $H^{-1}_{p}(\bR)$-valued function 
on $[0,1]$ we can write
$$
\frac{d}{dt}
\int_{\bR}u(t,x)\Phi_{x}(t,x)\,dx=
\frac{d}{dt}\langle u(t,\cdot),\Phi_{x}(t,\cdot)\rangle
=\langle u_{t}(t,\cdot),\Phi_{x}(t,\cdot)\rangle
$$
$$
+\langle u(t,\cdot),\Phi_{tx}(t,\cdot)\rangle,
$$
where $\langle\cdot,\cdot\rangle$ denotes the pairing between
$H^{1}_{q}(\bR)$ and $H^{-1}_{p}(\bR)$. Next we use the fact that
$\Phi_{tx}(t,\cdot)$ and $u_{x}(t,\cdot)$ are usual functions
(for almost all $t$), so that
$$
\langle u (t,\cdot),\Phi_{tx}(t,\cdot)\rangle
=\int_{\bR}u (t,x) \Phi_{tx}(t,x)\,dx
=-\int_{\bR}u_{x} (t,x) \Phi_{t}(t,x)\,dx.
$$
We also note that (for almost all $t$)
$$
\langle u_{t}(t,\cdot),\Phi_{x}(t,\cdot)\rangle
=\int_{\bR}f(t,x)\Lambda\Phi_{x}(t,x)\,dx
+\langle (au _{x}(t,\cdot))_{x},\Phi_{x}(t,\cdot)\rangle,
$$
where the last term equals
$$
-\langle u _{x}(t,\cdot) ,a\Phi_{xx}(t,\cdot)\rangle
= \int_{\bR}u _{x}(t,x)\Phi_{t}(t,x)\,dx.
$$

By combining these arguments we get that
$$
\frac{d}{dt}
\int_{\bR}u(t,x)\Phi_{x}(t,x)\,dx=
\int_{\bR}f(t,x)\Lambda\Phi_{x}(t,x)\,dx
$$
and for any $t_{0}\in[0,1)$
$$
\int_{\bR}u(t_{0},x)\Phi_{x}(t_{0},x)\,dx 
=\int_{0}^{t_{0}}\int_{\bR}f\Lambda\Phi_{x} \,dxdt.
$$
Finally, recall that $u=\Lambda w$ where $w\in\WO^{1,2}_{p}(Q)$,
which by Lemma 2.3.3 of \cite{LSU} implies that
$u$ is H\"older continuous in $x$ with the same exponent 
as in the proof of Theorem \ref{theorem 3.10.1}
and this allows us to get a contradiction by letting
$t_{0}\uparrow 1$ in the same way as in that proof. Both
theorems are thus proved.

 \mysection{Proof of Theorem \protect\ref{theorem 3.29.2}}
                                          \label{section 3.31.2}
Here we suppose that $\gamma\in(1,2)$, so that
 $\lambda\in(0,1)$.
One of solutions of \eqref{1.20.1} is
$$
\phi(x)=\phi_{\lambda}(x)=\int_{0}^{\infty}[1-e^{-2xr-r^{2}}]\frac{1}{r^{1+\lambda}}
\,dr ,
$$
with
$$
\phi'(x)=2\int_{0}^{\infty} e^{-2xr-r^{2}} \frac{1}{r^{ \lambda}}
\,dr>0,
$$
$$
\phi''(x)=-4\int_{0}^{\infty} e^{-2xr-r^{2}}  r^{1- \lambda}
\,dr<0.
$$
The fact that it is indeed a solution is seen from the following:
$$
\phi(x)=-\frac{1}{\lambda}\int_{0}^{\infty}
[1-e^{-2xr-r^{2}}] 
\,dr^{-\lambda}
$$
\begin{equation}
                                                          \label{5.7.5}
=\frac{1}{\lambda}\int_{0}^{\infty}(2x+2r)
e^{-2xr-r^{2}}r^{-\lambda}\,dr.
\end{equation}

Observe that, as is easily seen
(after substituting  $r=xs$), as $x\to\infty$
\begin{equation}
                                                \label{1.2.02}
\phi(x)\sim x^{\lambda}
\int_{0}^{\infty}[1-e^{-2 r }]\frac{1}{r^{1+\lambda}}
\,dr,\quad\phi'(x)\sim
2x^{\lambda-1}\int_{0}^{\infty} e^{-2 r } \frac{1}{r^{ \lambda}}
\,dr.
\end{equation}

Obviously $\phi(x)\to\pm\infty$ as $x\to\pm\infty$ and
$\phi(0)>0$. Therefore we can define
$x_{0}=x_{0\lambda}<0$ as a unique root of
$$
\phi(x_{0})=0
$$
and for $c>x_{0}$ let
$$
\psi_{c}(x)=\phi(x)\int_{c}^{x}\frac{1}{\phi^{2}(t)}
       e^{t^{2}}\,dt
$$
As in Section \ref{section 3.31.1} we   use this function
to investigate the behavior of $\phi(-x)$ as $x\to\infty$.

By Lemma \ref{lemma 1.3.1} the function
 $\psi_{c}(x)$ satisfies \eqref{1.20.1}
for any $c>c_{0}$. 
Also both functions $\psi_{-c_{0}}$ and $\phi(-x)$ vanish
at $x=c_{0}$ and hence there is a constant $m>0$ such that
$$
\phi(-x)=-m\psi_{-c_{0}}(x).
$$
In the same way as Lemma \ref{lemma 1.2.2}
is proved one gets the following.

\begin{lemma}
                                            \label{lemma 1.2.02}
We have
$$
\psi_{c}(x)\sim N\frac{1}{x^{1+\lambda}}e^{x^{2}} ,
\quad \phi(-x)\sim -N\frac{1}{x^{1+\lambda}}e^{x^{2}}
$$
as $x\to\infty$, where  the constants
 $N>0$   depend  only on  $\lambda$.

\end{lemma}

Next we set
$$
u(x)=\phi(x)+\phi(-x)
$$
 and in the same way as 
in Lemma
\ref{lemma 1.20.8} we prove that there exists 
a unique $c_{1}=c_{1\lambda}>0$
such that $(\phi (x)e^{-x^{2}})''=0$ at $x=c_{1}$.
In addition, 
\begin{equation}
                                                 \label{4.6.1}
(\phi (x)e^{-x^{2}})''<0\quad\text{for}\quad
0<x
< c_{1},\quad (\phi (x)e^{-x^{2}})''>0\quad\text{for}
\quad x> c_{1}.
\end{equation}

Furthermore as 
in Lemma
\ref{lemma 1.20.8}, the graph of the function $u (x)e^{-x^{2}}$
has inflection points on $(0,\infty)$ and we denote
by $c_{0\lambda}$ the smallest one.
Observe that 
\begin{equation}
                                                       \label{5.7.1}
u (c_{0\lambda})>0.
\end{equation}
 Indeed,
the equality $u (c_{0\lambda})=0$ is impossible due
to Remark \ref{remark 1.2.1}. However, if
$u (c_{0\lambda})<0$, then by Lemma \ref{lemma 1.2.1}
(ii) and Remark \ref{remark 1.2.1} the second-order derivative
of $u (x)e^{-x^{2}}$ at the closest zero of
 $u (x)$ lying to the left of $c_{0\lambda}$
is strictly positive and being negative at the origin
it would have another root smaller than $c_{0\lambda}$,
which contradicts its definition. Thus,
\begin{equation}
                                                 \label{4.6.2}
u (c_{0\lambda})>0,\quad
(u(x)e^{-x^{2}})''<0\quad\text{for}\quad
0<x
< c_{0\lambda}.
\end{equation}

By the way, notice that,   according to
 \eqref{1.2.02} and Lemma \ref{lemma 1.2.02} the function
$u_{\lambda}(x)e^{-x^{2}}$ approaches zero from
the negative side as $x\to\infty$. Therefore, 
there exists the smallest root $y>0$ of the equation
$u_{\lambda}(x)=0$ and as follows from the above $y>c_{0\lambda}$
and there are no inflection points between $c_{0\lambda}$ and $y$.
Then in the same way in which \eqref{1.22.1} is obtained
one shows that for $x\geq y$
$$
u_{\lambda}(x)=-\kappa\phi_{\lambda}(x)
\int_{y}^{x}\frac{1}{\phi_{\lambda}
^{2}(t)}e^{t^{2}}\,dt,
$$
where $\kappa\in(0,\infty)$ is a constant. It follows that
$u(x)<0$ for all $x>y$. Then the existence of a unique
root of the equation $(u_{\lambda}(x)e^{-x^{2}})''=0$
lying on $(y,\infty)$ is obtains as above, so that
$c_{0\lambda}$ is the smallest of the two roots.

Then we introduce $a_{1}=a_{1\gamma}$, $\beta(\gamma)=
a_{1\gamma}^{2}$,
 $w=w_{\gamma}$,
$\Psi=\Psi^{(\gamma)}$ and $a^{*}=a^{*}_{\gamma}$
in the same way as in \eqref{3.22.3}.

As before, \eqref{1.15.2} and \eqref{1.15.3} hold.
 However, in contrast with 
Lemma \ref{lemma 1.15.1} we now have the following.

\begin{lemma}
                                \label{lemma 1.15.01}
We have $c_{1}>c_{0}$, $a_{1}<1$, and
\begin{equation}
                                       \label{1.15.11}
\Psi_{t}(t,x)=\min_{a\in[a^{2}_{1} ,1]}[a\Psi_{xx}(t,x)].
\end{equation}

\end{lemma}

Proof. Assume that $c_{0}\geq c_{1}$. Then $a_{1}\geq1$
and  \eqref{1.15.2}, \eqref{1.15.3} imply that
$$
\Psi_{t}(t,x)=\max_{a\in[1, a^{2} ]}[a\Psi_{xx}(t,x)],
$$
$$
\Psi_{t}(t,x)\geq \Psi_{xx}(t,x).
$$
It follows by the maximum principle that
$$
\frac{1}{2\sqrt{\pi}}\int_{\bR}\Psi (t,y)
e^{-(x-y)^{2}/4}\,dy\leq \Psi(t+1,x),
$$
$$
2t^{-\gamma/2}
\int_{\bR}w(y/(2\sqrt{t}))e^{-y^{2}/4}\,dy\leq 2\sqrt{\pi}\Psi(t+1,0).
$$
However, the integral on the left is equivalent to
$$
4t^{(1-\gamma)/2}\int_{\bR}w(y)\,dy\to\infty
$$
as $t\downarrow0$. This yields a
contradiction, hence $c_{1}>c_{0}$, and the rest is trivial.
 The lemma is proved.

We are now ready to prove part of Theorem \ref{theorem 3.29.2}.
 \begin{theorem}
                              
There exists a constant   $\nu\in(0,1)$
 depending only on $\gamma\in(1,2)$
such that for any $\varepsilon\in(0,1/2)$
\begin{equation}
                                      \label{3.31.4}
\frac{1}{\varepsilon}\inf_{a\in \frA(a^{2}_{1\gamma},1)}
P_{a}
(1,0,0,[-\varepsilon,\varepsilon])
\geq \nu\frac{\varepsilon^{\gamma-1}}
{|\ln\varepsilon|^{\gamma/2}}.
\end{equation}

\end{theorem}

Proof. Take an $a\in \frA(a^{2}_{1\gamma},1)$, Since
$$
\Psi_{t}\leq a\Psi_{xx},
$$
we have 
$$
\int_{\bR}\Psi(t,y)\,P_{a}(1,x,0,dy)\geq \Psi(1+t,x),
$$
$$
\int_{\bR}\Psi(t,y)\,P_{a}(1,0,0,dy)\geq \Psi(1+t,0)\geq \Psi(2,0),
$$
where the last inequality holds for $t\in(0,1]$.
Here
$$
\int_{\bR}\Psi(t,y)\,P_{a}(1,0,0,dy)\leq
Mt^{-\gamma/2}P_{a}(1,0,0,[-2\sqrt{t|\ln t|},2\sqrt{t
|\ln t|}])
$$
$$
+t^{-\gamma/2}w(\sqrt{|\ln t|} ),
$$
where
$M=\max w$.
As $t\downarrow0$ we have
$$
t^{-\gamma/2}w(\sqrt{|\ln t|} )\leq N
t^{-\gamma/2}|\ln t|^{\lambda/2}e^{-|\ln t|}
=Nt^{1-\gamma/2}|\ln t|^{\lambda/2}\to0.
$$
It follows that there exists $t_{0}\in(0,1/2)$ such that
$$
t^{-\gamma/2}w(\sqrt{|\ln t|})\leq (1/2) \Psi(2,0)
$$
for $t\in(0,t_{0}]$. In that case upon setting
$\varepsilon=2\sqrt{t|\ln t|}$ we obtain
$$
\frac{1}{\varepsilon}
P_{a}(1,0,0,[-\varepsilon,\varepsilon])
\geq M^{-1}\Psi(2,0)t^{(\gamma-1)/2}
\frac{1}{2\sqrt{|\ln t|}}
$$
$$
=2^{-1}M^{-1}\Psi(2,0)\varepsilon^{\gamma-1}
\frac{1}{|\ln t|^{\gamma/2}}.
$$
To get \eqref{3.31.4}, now it only remains to observe that
$$
|\ln t|=2\ln2+2|\ln \varepsilon|+\ln|\ln t|\leq
2|\ln \varepsilon|+(1/2)|\ln t|
$$
for $t\in(0,t_{1}]$ with sufficiently small $t_{1}>0$,
so that 
$$
\frac{1}{|\ln t|^{\gamma/2}}\geq
\delta\frac{1}{|\ln \varepsilon|^{\gamma/2}}.
$$
The theorem is proved.

Now we are concerned with
another part  of Theorem  \ref{theorem 3.29.2}.
\begin{theorem}
                                  \label{theorem 1.13.1}
For  $\varepsilon\in(0,1 )$ 
set $a^{(\varepsilon)}(t,x)=a^{*}(\varepsilon+t,x)$.
Then 
we have
$$
\frac{1}{\varepsilon} 
P_{a^{(\varepsilon)}}(1,0,0,[-\varepsilon,\varepsilon])\leq \mu
^{-1} \Psi(1,0)
\varepsilon^{\gamma-1},
$$
where $\mu=\inf_{|x|\leq 1/2} w(x)$ depends only on $\gamma$.
\end{theorem}

Proof. Indeed, we have
$$
\int_{\bR}\Psi(t,y)\,P_{a^{(\varepsilon)}}(1,0,0,dy)= \Psi(1+t,0)\leq
\Psi(1,0),
$$
$$
t^{-\gamma/2}
\mu P_{a^{(\varepsilon)}}(1,0,0, [-\sqrt{t},\sqrt{t}])\leq \Psi(1,0).
$$

This obviously leads to the desired result
and proves the theorem.

Similarly to Corollary \ref{corollary 4.19.1}
we have
\begin{corollary}
                                            \label{corollary 4.19.2}
The function $\beta(\gamma)=a^{2}_{1\gamma}$
is a  decreasing function of $\gamma\in(1,2)$.
\end{corollary}

To finish the proof of Theorem  \ref{theorem 3.29.2}
we only need to show that
  $\beta(\gamma)=a^{2}_{1\gamma}\to1$ 
as $\gamma\downarrow1$
and $\beta(\gamma)$ tends to a nonzero limit as $\gamma
\uparrow 2$. We already know that these limits exist.

Observe that
the condition $(ve^{-x^{2}})''=0$
for solutions $v>0$ of \eqref{1.20.1} 
 is also written as
$$
\frac{v''(x)}{2v(x)}=2x^{2}-1-2\lambda.
$$
Therefore, $c_{0}=c_{0\gamma}$ and  $c_{1}=c_{1\gamma}$
also satisfy
\begin{equation}
                                                        \label{5.7.2}
\frac{\phi_{\lambda}''(c_{0})+\phi_{\lambda}''(-c_{0})}
{ 2\phi_{\lambda}(c_{0})+2\phi_{\lambda}(-c_{0})}
=2c_{0}^{2}-1-2\lambda,\quad
\frac{\phi_{\lambda}''(c_{1}) }
{ 2\phi_{\lambda}(c_{1}) }
=2c_{1}^{2}-1-2\lambda.
\end{equation}
Here  $\phi_{\lambda}''\leq0$ and $\phi_{\lambda}(c_{0})
+ \phi_{\lambda}(-c_{0})>0$
owing to \eqref{5.7.1}.
Hence $c_{0\gamma}$ and $c_{1\gamma}$
are bounded. Also   $\phi_{\lambda}''$ is bounded on bounded
intervals and $\phi_{\lambda}\to\infty$ uniformly on  intervals
of type $[0,b]$, $b>0$,
as $\gamma=\lambda+1\downarrow1$ because of the divergence
at infinity
of the integral defining $\phi_{\lambda}$. This immediately
implies that
\begin{equation}
                                                      \label{5.7.3}
\lim_{\gamma\downarrow1}c_{1\gamma}=2^{-1/2}
\end{equation}
and, along with Corollary \ref{corollary 4.19.2}, shows 
 that the limit of $c_{0\gamma}$ as $\gamma\downarrow1$
exists. We denote it by $c$. By using 
\eqref{5.7.5} we rewrite the first equation in \eqref{5.7.2} as
$$
\lambda[\phi_{\lambda}''(c_{0})+\phi_{\lambda}''(-c_{0})]
$$
$$
=2[2c_{0}^{2}-1-2\lambda]
\int_{0}^{\infty}[(2c_{0}+2r)e^{-2c_{0}r}+(2r-c_{0})
e^{2c_{0}r}]e^{-r^{2}}r^{-\lambda}\,dr.
$$
By letting $\lambda\downarrow0$ we obtain
\begin{equation}
                                                      \label{5.7.6}
2[2c ^{2}-1 ]
\int_{0}^{\infty}[(2c +2r)e^{-2c r}+(2r-c )
e^{2c r}]e^{-r^{2}} \,dr=0.
\end{equation}
As a function, say $f$, of $c$ the last integral is obtained as the limit
as $\lambda\downarrow0$ of solutions of the equation
$v''-2xv'+\lambda v=0$. Therefore, $f''-2xf'=0$,
$$
f'=C_{1}e^{ x^{2}},\quad f(x)=C_{1}\int_{0}^{x}e^{ t^{2}}\,dt
+C_{2},
$$
where $C_{1},C_{2}$ are some constants. The function $f$
is obviously even, hence $C_{1}=0$. Also $C_{2}=f(0)>0$
and \eqref{5.7.6} implies that 
$$
c=\lim_{\gamma\downarrow1}c_{0\gamma}=2^{-1/2},\quad
\lim_{\gamma\downarrow1}\beta(\gamma)=1.
$$

Next, by comparing \eqref{5.7.2} with \eqref{4.4.2},
 the left-hand sides of which
is positive in our situation, we obtain
$$
0\leq 2c_{1}^{2}-1-\lambda\leq\lambda,\quad
1+\lambda<2c_{1}^{2}<1+2\lambda
$$
which implies \eqref{5.7.3} one more time and also shows that
   $c_{1\gamma}$ is separated away from zero.
 On the set of possible values of 
$c_{1\gamma}$ we have that $\phi_{\lambda}\to\infty$  because of 
the divergence at zero
of the above mentioned integral as $\lambda\uparrow1$. Hence,
by \eqref{5.7.2}
  $$
\lim_{\gamma\uparrow2}c_{1\gamma}=(3/2)^{1/2}.
$$

Again this implies  
that the limit of $c_{0\gamma}$ as $\gamma\uparrow2$
exists. Moreover,
as is easy to see,
$\phi_{\lambda}(x)+\phi_\lambda(-x)$ tends to a finite
limit, say $\psi$, as $\lambda\uparrow1$ and therefore $c_{0\gamma}$
tends to, say, $d\geq0$ satisfying
$$
 \psi ''(d )   
= 2\psi ( d)(2d ^{2}-3).
$$ 
Since $\phi_{\lambda}(x)+\phi_\lambda(-x)$ are solutions 
of \eqref{1.20.1} we have that $\psi''-2x\psi'+2\psi=0$,
which shows that
$\psi''(0)=-2 \psi (0)$.
In particular, $d\ne0$
  and this finishes the proof of the theorem.

\end{document}